\documentclass[12pt]{amsart}
\usepackage{enumerate}
\usepackage{latexsym}
\usepackage[pdftex]{graphicx}
\usepackage{amssymb}
\usepackage{comment}
\usepackage{xcolor}

\newtheorem*{theorem}{Theorem}
\newtheorem{lemma}{Lemma}
\newtheorem{proposition}{Proposition}
\newtheorem{remark}{Remark}
\newtheorem{example}{Example}

\newtheorem{definition}{Definition}
\newtheorem{corollary}{Corollary}

\newtheorem{problem}{Problem}

\newcommand{\ve}{\varepsilon}

\def\R{{\mathbb R}}

\newcommand{\beq}{\begin{equation}}
\newcommand{\eeq}{\end{equation}}
\newcommand{\beqna}{\begin{eqnarray*}}
\newcommand{\eeqna}{\end{eqnarray*}}
\newcommand{\beqn}{\begin{equation*}}
\newcommand{\eeqn}{\end{equation*}}
\newcommand{\bp}{\begin{proof}}
\newcommand{\ep}{\end{proof}}
\newcommand{\bprop}{\begin{proposition}}
\newcommand{\eprop}{\end{proposition}}
\newcommand{\bt}{\begin{theorem}}
\newcommand{\et}{\end{theorem}}
\newcommand{\bex}{\begin{example}}
\newcommand{\eex}{\end{example}}
\newcommand{\bc}{\begin{corollary}}
\newcommand{\ec}{\end{corollary}}
\newcommand{\bl}{\begin{lemma}}
\newcommand{\el}{\end{lemma}}
\newcommand{\bprob}{\begin{problem}}
\newcommand{\eprob}{\end{problem}}
\newcommand{\br}{\begin{remark}}
\newcommand{\er}{\end{remark}}
\newcommand{\bd}{\begin{definition}}
\newcommand{\ed}{\end{definition}}

\renewcommand{\ge}{\geqslant}
\renewcommand{\le}{\leqslant}

\newcommand{\KR}{\mathcal K_R}

\begin{document}

\title
[Fine approximation of  convex bodies by polytopes]
{Fine approximation of  convex bodies by polytopes}

\author[M. Nasz\'odi]{M\'arton Nasz\'odi}
\address{Dept. of Geometry, Lorand E\"otv\"os University\\
Pazmany Peter Stny. 1/C, Budapest, Hungary 1117} 
\email{marton.naszodi@math.elte.hu}

\author[F. Nazarov]{Fedor Nazarov}
\address{Department of Mathematics, Kent State University,
Kent, OH 44242, USA} \email{nazarov@math.kent.edu}

\author[D. Ryabogin]{Dmitry Ryabogin}
\address{Department of Mathematics, Kent State University,
Kent, OH 44242, USA} \email{ryabogin@math.kent.edu}

\thanks{The first named author is  
supported in part by the J\'anos Bolyai Research Scholarship of 
the Hungarian Academy of Sciences, and the National Research, Development, and 
Innovation Office, NKFIH Grant K119670.}

\thanks{The second and the third named authors are  
supported in part by U.S.~National Science Foundation Grants  DMS-0800243 and 
DMS-1600753}


\begin{abstract}
We prove that for every convex body $K$ with the center of mass at the origin 
and every  $\varepsilon\in \left(0,\frac{1}{2}\right)$, 
there exists a convex polytope $P$ with at most 
$e^{O(d)}\varepsilon^{-\frac{d-1}{2}}$ vertices such that 
$(1-\varepsilon)K\subset P\subset K$.
\end{abstract}

\maketitle

\section{Introduction and main result}

A \emph{convex body} in 
${\mathbb R^d}$ is a compact convex set with non-empty interior.
Our goal is to prove the following theorem.

\bt\label{FM}
Let $K$ be a convex body in ${\mathbb R^d}$ with the center of mass at the 
origin, and let $\varepsilon\in \left(0,\frac{1}{2}\right)$. 
Then there exists a convex polytope $P$ with at most 
$e^{O(d)}\varepsilon^{-\frac{d-1}{2}}$ vertices such that 
$(1-\varepsilon)K\subset P\subset K$.
\et

This result improves the 2012 theorem of Barvinok \cite{B} by removing the 
symmetry assumption and the extraneous $(\log\frac{1}{\varepsilon})^d$ factor. 
Our approach uses a mixture of geometric and probabilistic tools.

We refer the reader to the surveys of Bronshtein \cite{Br} and Gruber \cite{Gr} for the discussion of the history of the 
problem. Unfortunately, we will have to rely upon two non-trivial classical 
results (Blaschke-Santal\'o inequality and its reverse), which makes this paper a bit less reader-friendly than we would like 
it to be despite our best efforts to provide well-written and easily accessible 
references for all statements that we use without a proof.

\section{Outline of the proof}

Without loss of generality, we may assume that $K$ has smooth boundary, in particular, $K$ has a unique supporting hyperplane at 
each boundary point. Our task is to find a finite set of points 
$Y\subset\partial K$
such that $P=\operatorname{conv} Y$ satisfies $(1-\ve)K\subset P$. By 
duality, this is equivalent to the requirement that 
every
cap $S(x,\ve)=\{y\in\partial K: \langle y,\nu_x\rangle\ge (1-\ve)\langle 
x,\nu_x\rangle\}$, where $x\in\partial K$ and $\nu_x$ is the outer unit normal to
$\partial K$ at $x$, contains at least one point of $Y$. 

The key idea is to construct a probability measure $\mu$ on $\partial K$
such that for every $x\in\partial K,\ve\in \left(0,\frac 12\right)$, we have 
$\mu(S(x,\ve))\ge p\ve^{\frac{d-1}2}$ with some $p=e^{O(d)}$ depending
on $d$ only. 

Since  there are infinitely many caps,
 our next aim is to choose an appropriate finite net $X\subset \partial K$ 
of cardinality $C(d)\ve^{-\frac{d-1}2}$
such that the condition $S(x,\frac \ve 2)\cap Y\ne \varnothing$ for all $x\in 
X$ implies that
$S(x,\ve)\cap Y\ne \varnothing$ for all $x\in \partial K$. 
Given such a net, we will be able to apply a general combinatorial result essentially due to Rogers to construct the desired  set $Y$ of cardinality approximately
$\log 
C(d)p^{-1}\ve^{-\frac{d-1}2}$, which will be still $e^{O(d)}\ve^{-\frac{d-1}2}$ as long as $C(d)$ is at most double 
exponential in $d$.

A natural net to try is the Bronshtein--Ivanov net (see 
\cite{BI}), which allows one to approximate a point $x\in\partial K$ 
and the corresponding outer unit normal $\nu_x$ by a point in the net and its 
outer unit normal simultaneously. Unfortunately,
it works only for uniformly $2$-convex bodies, i.e., the bodies that can be 
touched by an outer sphere of fixed controllable radius at every boundary 
point.

So, the last step will be to show that the task of approximating an 
arbitrary convex body $K$ can be reduced to that of approximating a certain 
uniformly $2$-convex body associated with $K$.

In the exposition, these steps are presented in reverse. We start with 
constructing the associated uniformly $2$-convex body
(Sections \ref{Sp}, \ref{Tf}, \ref{Ft}). Then we build the Bronshtein-Ivanov 
net $X$ of appropriate mesh and cardinality,
and check that it is, indeed, enough to consider the caps $S(x,\frac \ve 2), 
x\in X$ (Sections \ref{TBI}, \ref{Td}, \ref{D}).
Finally, we construct the probability measure $\mu$ and finish the proof of the theorem
(Sections \ref{TR}, \ref{Tc}).

\section{Standard position}
\label{Sp}

Since the problem is invariant under linear transformations, we can always 
assume that our body $K$ is in some ``standard position". The exact notion of 
the standard position  to use is not very important as long as it guarantees 
that $B\subset K\subset d^2 B$, say, where $B$ is the unit ball in ${\mathbb 
R^d}$ centered at the origin.

One possibility is to make a linear transformation such that John's ellipsoid 
(see \cite{Ba}, Lecture 3) of the centrally-symmetric convex body $L=K\cap\,-K$ 
is the unit ball, so $B\subset L\subset\sqrt{d}B$ and, since 
$K\supset\,-\frac{1}{d}K$ (see \cite{BF}, page 57), it follows that $B\subset 
K\subset d\sqrt{d} B$.

\section{The  function $\varphi_{\delta}$ and the mapping $\Phi_{\delta}$}
\label{Tf}

Fix $\delta\in \left(0,\frac{1}{2}\right)$. For $r\ge 0$, define 
$\varphi_{\delta}(r)$ as the positive root of the equation $\varphi+\delta 
r^2\varphi^2=1$. Put $\Phi_{\delta}(x)=x\varphi_{\delta}(|x|)$, $x\in {\mathbb 
R^d}$.

\begin{figure}[ht]
\includegraphics[height=3.5cm]{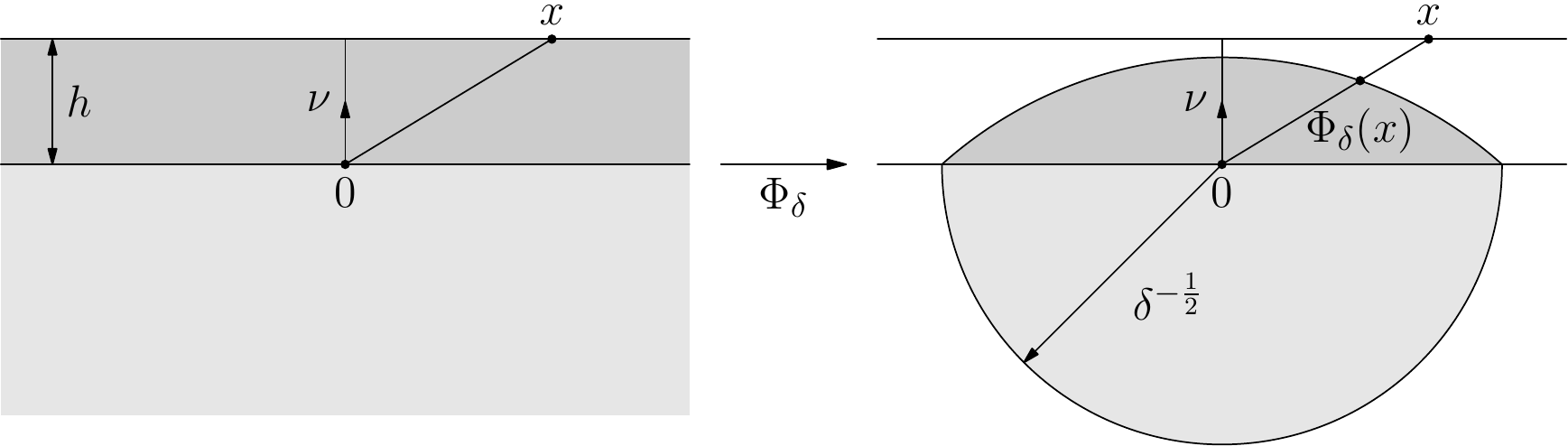}
\caption{The mapping $\Phi_\delta$}
\label{fig1}
\end{figure}

\bl\label{fatn}
The function $\varphi_{\delta}$ is a decreasing smooth function on 
$[0,+\infty)$; $r\mapsto r\varphi_{\delta}(r)$ is an increasing function 
mapping  $[0,+\infty)$ to $[0,\delta^{-\frac{1}{2}})$; $\Phi_{\delta}$ is a 
diffeomorphism of ${\mathbb R^d}$ onto the open ball $\delta^{-\frac{1}{2}} 
\operatorname{int} B$; if $\nu$ is a unit vector and $h>0$, then the image 
$\Phi_{\delta}(H_{\nu,h})$ of the half-space $H_{\nu,h}=\{x:\langle 
x,\nu\rangle\le h\}$ is the intersection of $\delta^{-\frac{1}{2}} 
\operatorname{int} B$ and the ball of radius $\sqrt{\frac{1}{4\delta^2 
h^2}+\frac{1}{\delta}}$ centered at $-\frac{1}{2\delta h}\nu$
(see Figure \ref{fig1}).\el

\bp
The first statement is obvious. To show the second one, just notice that 
$r\varphi_{\delta}(r)$ is the positive root of $\frac{\psi}{r}+\delta\psi^2=1$ 
and, as $r\to\infty$, this root increases to $\delta^{-\frac{1}{2}}$. The third 
claim follows from the observation that the derivative of the mapping 
$r\mapsto r\varphi_{\delta}(r)$ is strictly positive and continuous on 
$[0,+\infty)$. 
To prove the last claim, observe that if $\langle x,\nu\rangle=h$, then
\begin{multline*}
\Big|\Phi_{\delta}(x)+\frac{1}{2\delta 
h}\nu\Big|^2=\Big|x\varphi_{\delta}(|x|)+\frac{1}{2\delta h}\nu\Big|^2=
\\
|x|^2\varphi_{\delta}(|x|)^2+\frac{\varphi_{\delta}(|x|)}{\delta}+\frac{1}{
4\delta^2 h^2}=\frac{1}{4\delta^2 h^2}+\frac{1}{\delta}
\end{multline*}
by the definition of $\varphi_{\delta}$. \ep

It follows that for every convex body $K$ containing the origin, 
$\Phi_{\delta}(K)$ is also convex. 
Since for every interval $I_x=\{rx:0\le r\le 1\}$, $x\in\mathbb R^d$, we have 
$\Phi_\delta(I_x)\subset I_x$, the image $\Phi_\delta(K)$ is contained in $K$.
Moreover, if $B\subset K$, then $\Phi_{\delta}(K)$ is the intersection of balls 
of radii not exceeding $\sqrt{\frac{1}{4\delta^2 
}+\frac{1}{\delta}}\le\frac{1}{\delta}$. In particular, for every boundary  
point $x\in \partial \Phi_{\delta}(K)$, we can find a ball of radius 
$\frac{1}{\delta}$ containing $\Phi_{\delta}(K)$ whose boundary sphere touches 
$\Phi_{\delta}(K)$ at $x$.

\section{From the approximation of $\Phi_{\delta}(K)$ to the approximation of 
$K$}
\label{Ft}

\bl\label{po}
Let $\varepsilon\in \left(0,\frac{1}{2}\right)$. Suppose that a convex body $K$ 
satisfies $0\in K\subset d^2 B$ and $\delta<\frac{1}{4d^4}$. If 
$Y\in\partial K$ is a finite set such that  $\left(1-\frac{\varepsilon}{2}\right)\Phi_{\delta}(K)
\subset \operatorname{conv}(\Phi_{\delta}(Y) )$, then 
$(1-\varepsilon)K\subset \operatorname{conv}(Y)$.
\el
\bp
Note that the conditions of the lemma imply that 
$0\in\operatorname{conv}
(\Phi_{\delta}(Y))$. Since for every $y\in{\mathbb R^d}$, 
$\Phi_\delta(y)$ is a positive multiple of $y$, we conclude that $0\in 
P=\operatorname{conv}(Y)$ as well,
so $\Phi_\delta(P)$ is convex. 
Suppose that there exists $x\in K$ such that $(1-\varepsilon)x\notin P$. Then,
$$
\Phi_{\delta}((1-\varepsilon)x)\notin \Phi_{\delta}(P)\supset 
\operatorname{conv}(\Phi_{\delta}(Y)).
$$
However, 
$$
\Phi_{\delta}((1-\varepsilon)x)=(1-\varepsilon)\frac{\varphi_{\delta}
((1-\varepsilon)|x|)}{\varphi_{\delta}(|x|)}\Phi_{\delta}(x).
$$
Denoting $\eta_t=\varphi_{\delta}((1-t)|x|)$, $t\in [0,1]$, we have 
$$
\eta_{\varepsilon}+\delta(1-\varepsilon)^2|x|^2\eta_{\varepsilon}
^2=\eta_0+\delta|x|^2\eta_0^2=1.
$$
Since $\,\delta |x|^2\eta_{\varepsilon}^2\ge \delta |x|^2\eta_0^2\,$ and 
$\,\delta \varepsilon^2|x|^2\eta_{\varepsilon}^2\ge 0$, it follows that
$$
\eta_{\varepsilon}(1-2\delta \varepsilon|x|^2\eta_{\varepsilon})\le 
\eta_0,\quad\textrm{so}\quad 
\frac{\eta_{\varepsilon}}{\eta_0}\le\frac{1}{1-2\delta 
\varepsilon|x|^2\eta_{\varepsilon}}.
$$
Since $\eta_{\varepsilon}\le1$ and $2\delta|x|^2\le 2\delta d^4\le 
\frac{1}{2}$, we get 
$$
(1-\varepsilon)\frac{\eta_{\varepsilon}}{\eta_0}\le 
\frac{1-\varepsilon}{1-\frac{\varepsilon}{2}}\le 1-\frac{\varepsilon}{2},
$$
 so  $\left(1-\frac{\varepsilon}{2}\right)\Phi_{\delta}(x)$ cannot be contained 
in  $\operatorname{conv}(\Phi_{\delta}(Y))$, which 
contradicts our assumption.
\ep

This lemma implies that an $\frac{\varepsilon}{2}$-approximation of $\Phi_{\delta}(K)$ yields an 
$\varepsilon$-approximation of $K$. Note also that  $\Phi_{\delta}(K)$  is rather close to $K$. More precisely, if $0\in K\subset d^2B$, we have 
$(1-\delta d^4)K\subset  \Phi_{\delta}(K)\subset K$. The center of mass of  $\Phi_{\delta}(K)$   may no longer be at the origin, of course, but the only non-trivial property of $K$  we shall really use is the Santal\'o bound  $\textrm{vol}_d(K)\textrm{vol}_d(K^{\circ})\le e^{O(d)}d^{-d}$, where 
$$
K^{\circ}=\{y\in{\mathbb R^n}:\langle x, y\rangle\le 1\quad\textrm{for 
all}\quad x\in K\}
$$
is the polar body of the convex body $K$. This bound holds for $K$ because $0$, being the center of mass of $K$, is, thereby, the 
Santal\'o point of $K^{\circ}$ (see Section \ref{Tc} for details). For sufficiently small $\delta>0$, it is inherited by 
$\Phi_{\delta}(K)$ just because $(\Phi_{\delta}(K))^{\circ}\subset (1-\delta d^4)^{-1}K^{\circ}$ and, thereby,
$$
\textrm{vol}_d(\Phi_{\delta}(K))\textrm{vol}_d((\Phi_{\delta}(K))^{\circ})\le (1-\delta d^4)^{-d}\textrm{vol}_d(K)\textrm{vol}_d(K^{\circ}).
$$
Choosing $\delta=\frac{1}{4d^5}$, we see that the body $\Phi_{\delta}(K)$ also satisfies the Santal\'o bound with only marginally worse constant. At last, if $B\subset K$, we have 
$\frac{1}{2}B\subset (1-\delta)B\subset \Phi_{\delta}(K)$. 

Thus, replacing $K$ by $\Phi_{\delta}(K)$ (and $\varepsilon$ by $\frac{\varepsilon}{2}$) if necessary, from now on we can restrict ourselves to the 
class ${\mathcal K}_R$ of convex bodies $K$ with smooth boundary such that $\frac{1}{2}B\subset K\subset d^2B$ and for every boundary point 
$x\in\partial K$, there exists a  ball of fixed radius $R=4d^5$ containing $K$ whose boundary sphere touches $K$ at $x$. Moreover, we can also assume that 
$\textrm{vol}_d(K)\textrm{vol}_d(K^{\circ})\le e^{O(d)}d^{-d}$.

\section{The Bronshtein--Ivanov net}
\label{TBI}

\begin{figure}[ht]
\includegraphics[height=4.5cm]{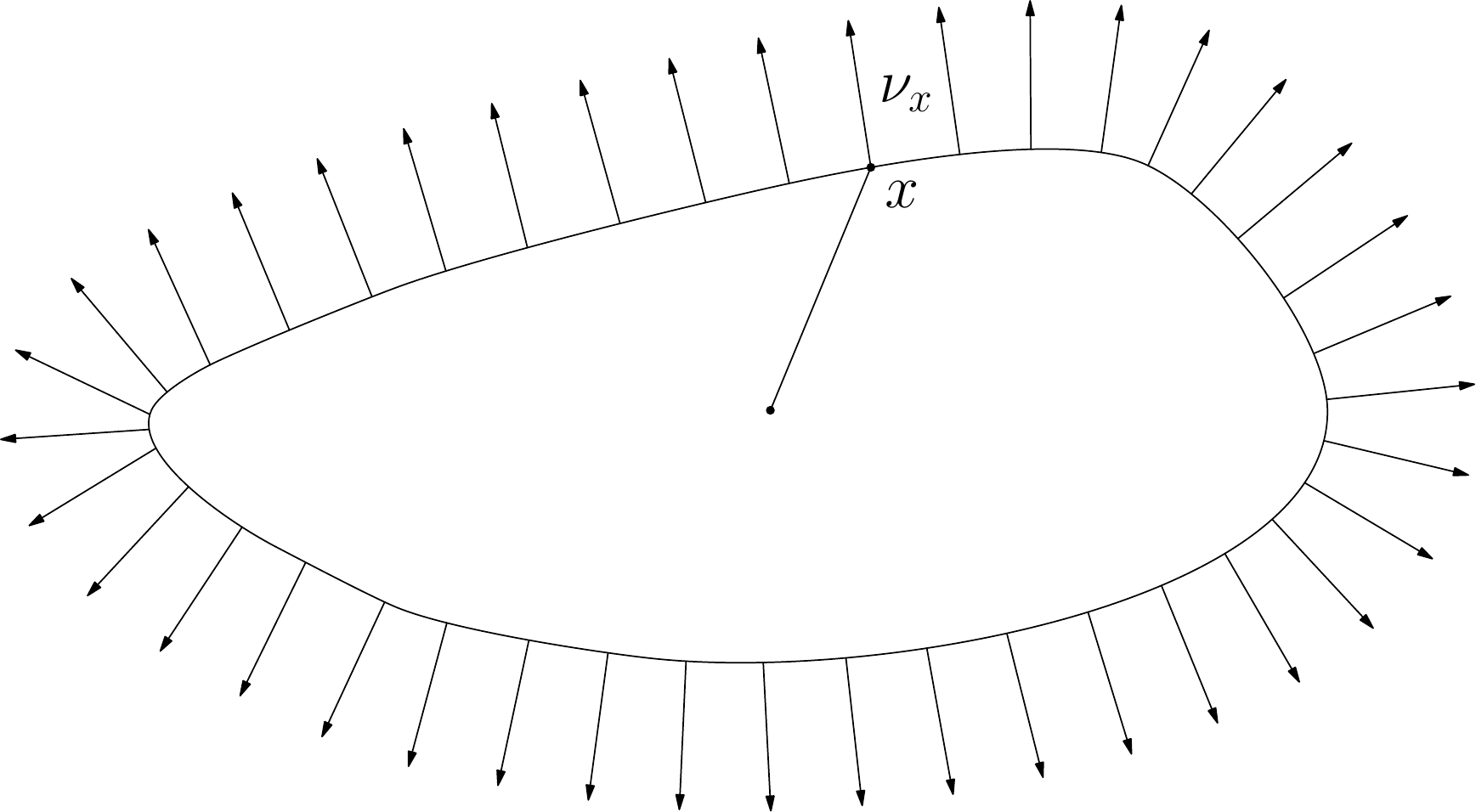}
\caption{The Bronshtein-Ivanov net}
\label{fig2}
\end{figure}

Let $\rho\in\left(0,\frac{1}{2}\right)$. Let $K$ be a convex body with   smooth 
boundary containing the origin and contained in $d^2B$. Consider the set $S$ of 
points $\{x+\nu_x: x\in \partial K\}$, where $\nu_x$ is the outer unit normal 
to $\partial K$ at $x$. Let $\{x_j+\nu_{x_j}:1\le j\le N\}$ be a maximal 
$\rho$-separated set in $S$, i.e., a set such that any two of 
its members are at distance at least $\rho$ (see Figure \ref{fig2}).
We will call the corresponding set $\{x_j:1\le j\le N\}$ a Bronshtein-Ivanov net of mesh $\rho$ for the body $K$.

\bl\label{bi}
We have $N\le 2^d(d^2+3)^d\rho^{-d+1}$, and for every $x\in\partial K$, we can 
find $j$ such that $|x-x_j|^2+|\nu_x-\nu_{x_j}|^2\le \rho^2$.
\el
\bp
Let $x', x''\in\partial K$ and let $\nu'=\nu_{x'}$, $\nu''=\nu_{x''}$. Note 
that, by the convexity of $K$, we must have $\langle\nu',x'-x''\rangle\ge 0$, 
$\langle\nu'',x''-x'\rangle\ge 0$.
Hence,  we always have
\begin{multline*}
|x'+\nu'-x''-\nu''|^2=
\\
|x'-x''|^2+|\nu'-\nu''|^2+2(\langle\nu',x'-x''\rangle+
\langle\nu'',x''-x'\rangle)\ge
\\
 |x'-x''|^2+|\nu'-\nu''|^2,
\end{multline*}
and the second conclusion  of the lemma follows immediately from the definition 
of $x_j$.

\begin{figure}[ht]
\includegraphics[height=5cm]{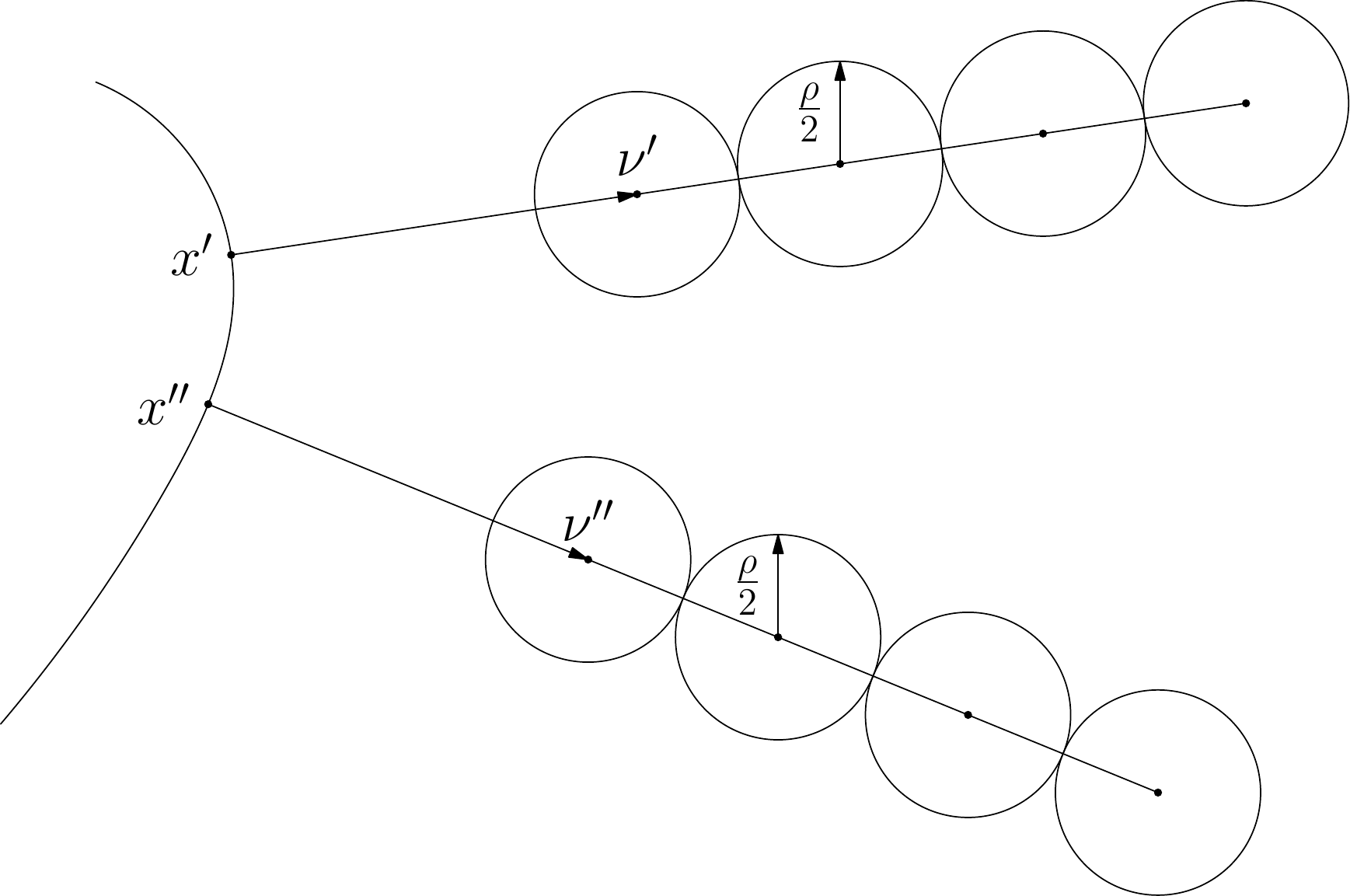}
\caption{The disjoint balls}
\label{fig3}
\end{figure}

Now assume that $s',s''\ge 0$. Write
\begin{multline*}
|x'+\nu'+s'\nu'-x''-\nu''-s''\nu''|^2=|x'+\nu'-x''-\nu''|^2+
\\
|s'\nu'-s''\nu''|^2+2s'\langle\nu',x'-x''\rangle+2s''\langle\nu'',x''-x'\rangle 
+
\\
2(s'+s'')(1-\langle\nu',\nu''\rangle)\ge
|x'+\nu'-x''-\nu''|^2.
\end{multline*}
Thus,  if the balls of radius $\frac{\rho}{2}$ centered at $x'+\nu'$ and 
$x''+\nu''$ are disjoint, so are the balls of radius $\frac{\rho}{2}$ centered 
at $x'+(1+s')\nu'$ and $x''+(1+s'')\nu''$. From here we conclude that the balls 
of radius $\frac{\rho}{2}$ centered at the points $x_j+(1+k\rho)\nu_{x_j}$, 
$0\le k\le \frac{1}{\rho}$ are all disjoint (see Figure \ref{fig3}) 
and contained in $(d^2+3)B$.

 The 
total
number of these balls is at least $\frac N\rho$ (for every point $x_j$ in the 
net, there is a chain of at least $\frac 1\rho$ balls 
corresponding to different values of $k$), whence 
$\frac{N}{\rho}\le\left(\frac{d^2+3}{\frac{\rho}{2}}\right)^d$ and the desired 
bound for $N$ follows.
\ep

\section{The distance bound}
\label{Td}

The following lemma shows that $\varepsilon$-caps of convex bodies $K\in {\mathcal K}_R$ have small diameters.

\bl\label{bpnx}
Let $\varepsilon\in\left(0,\frac 12\right)$. Assume that $K\in \KR$, $x\in\partial K$, and $\nu$ is the outer normal  to $\partial K$ at $x$. If $y\in S(x,\varepsilon)$, i.e., $y\in K$ and 
 $\langle y,\nu\rangle\ge (1-\varepsilon)\langle x,\nu\rangle$, then
 $|y-x|\le 
\sqrt{2R}\,d\,\sqrt{\varepsilon}$.
\el
\bp
Let $Q$ be the ball of radius $R$ containing $K$ whose boundary sphere touches 
$K$ at $x$. Then $y\in Q$ and $\nu$ is the outer unit normal to $Q$ at $x$, so 
$Q$ is centered at $x-R\nu$.
Note also that, since $0\in K\subset d^2 B$, we have $0\le \langle 
x,\nu\rangle\le  d^2$. Now we have
$$
R^2\ge |y-x+R\nu|^2=|y-x|^2+2R\langle y-x,\nu\rangle+R^2,
$$
so
$$
|y-x|^2\le 2R\langle x-y,\nu\rangle\le 2R\varepsilon \langle x,\nu\rangle\le 2R 
d^2\varepsilon,
$$
as required.
\ep

\section{Discretization}
\label{D}

\bl\label{dv}
Let $\varepsilon,\rho\in\left(0,\frac{1}{2}\right)$. Let $K\in \KR$. Let $x, 
x', y\in\partial K$ and let $\nu$ and $\nu'$ be the outer unit normals to 
$\partial K$ at $x$ and $x'$ respectively. Assume that 
$|x-x'|^2+|\nu-\nu'|^2\le \rho^2$ and
$\langle y,\nu\rangle\ge \left(1-\frac{\varepsilon}{2}\right)\langle 
x,\nu\rangle$. Then 
$$
\langle y,\nu'\rangle\ge \left(1-\frac{\varepsilon}{2}-2\rho(\rho+\varepsilon 
d^2+|y-x|)\right)\langle x',\nu'\rangle.
$$
\el
\bp
We have
\begin{multline*}
\langle y,\nu'\rangle=\langle x,\nu'\rangle+\langle y-x,\nu'\rangle=
\\
\langle x',\nu'\rangle+\langle x-x',\nu'\rangle+\langle y-x,\nu\rangle+\langle 
y-x,\nu'-\nu\rangle\ge
\\
\langle x',\nu'\rangle+\langle x-x',\nu'-\nu\rangle+\langle 
y-x,\nu\rangle+\langle y-x,\nu'-\nu\rangle\ge
\\
\langle x',\nu'\rangle-\rho^2-\frac{\varepsilon}{2}\langle 
x,\nu\rangle-\rho|y-x|.
\end{multline*}
Here, when passing from the second line to the third one,
we used the inequality $\langle x-x', \nu\rangle\ge 0$.

Note now that
$$
\langle x,\nu\rangle=\langle x,\nu'\rangle+\langle x,\nu-\nu'\rangle\le \langle 
x',\nu'\rangle+\rho d^2
$$
and $\langle x',\nu'\rangle\ge\frac{1}{2}$, so 
\begin{multline*}
\langle y,\nu'\rangle\ge \Big(1-\frac{\varepsilon}{2}\Big)\langle 
x',\nu'\rangle-\rho\Big(\rho+\frac{\varepsilon d^2}{2}+|y-x|\Big)\ge 
\\
\Big(1-\frac{\varepsilon}{2}-2\rho(\rho+\varepsilon d^2+|y-x|)\Big)\langle 
x',\nu'\rangle.
\end{multline*}
\ep

Recall that our task is to find a finite set of points $Y\subset\partial K$ 
such that $(1-\varepsilon)K
\subset \operatorname{conv} Y$. This requirement is equivalent to the statement 
that for every 
$x\in \partial K$, there exists $y\in Y$ such that $\langle y,\nu\rangle\ge 
(1-\varepsilon) \langle x,\nu\rangle$, where $\nu$ is the outer unit normal to 
$\partial K$ at $x$.

Lemma \ref{dv} implies that it would suffice  to show the existence of $y\in Y$ 
satisfying a slightly stronger inequality 
$\langle y,\nu\rangle\ge  \left(1-\frac{\varepsilon}{2}\right) \langle 
x,\nu\rangle$ for every
 point $x$ in the Bronshtein--Ivanov net only, provided that  we can ensure 
that $2\rho(\rho+\varepsilon d^2+|y-x|)\le\frac{\varepsilon}{2}$.

To this end, we apply Lemma  \ref{bpnx}, which shows that  the inequality 
$\langle y,\nu\rangle\ge  \left(1-\frac{\varepsilon}{2}\right) \langle 
x,\nu\rangle$ automatically implies the distance bound
$|y-x|\le \sqrt{2R}d\sqrt{\frac{\varepsilon}{2}}= 
d\sqrt{R}\,\sqrt{\varepsilon}$.
Thus, if we choose $\rho=\frac{1}{4(d^2+1+d\sqrt{R})}\sqrt{\varepsilon}$, we 
will be in good shape.

By Lemma~\ref{bi}, the  size $N$ of the corresponding 
Bronshtein-Ivanov net  is at most
$8^d(d^2+3)^d(d^2+1+d\sqrt{R})^d\varepsilon^{-\frac{d-1}{2}}=C(d)\varepsilon^{
-\frac{d-1}{2}}$, which has the correct power of $\varepsilon$ already. 
However, 
$C(d)$ is superexponential in $d$, which prevents us from just using the full 
Bronshtein--Ivanov net for $Y$ and forces us to work a bit harder.

\section{Rogers' trick}
\label{TR}

We now remind the reader a simple abstract construction essentially due to Rogers \cite{R}.

\bl\label{lem:transversal}
Let ${\mathcal S}=\{S_1,\dots,S_N\}$ be a family of measurable subsets of a probability space $(U,\mu)$ such that for some $\vartheta>0$,  we have $\mu(S_i)\ge \vartheta$ for all $i=1,\dots,N$. Then there exists a set $Y$ of cardinality at most 
$\lceil \vartheta^{-1}\log (N\vartheta)\rceil+\vartheta^{-1}$ that intersects each $S_i$.
\el
Here $\lceil z\rceil$ stands for the least non-negative integer greater than or equal to $z$.
\bp
First we choose $M$ points randomly and independently according to $\mu$ and obtain a random set $Y_0$. For every fixed $i\in\{1,\dots, N\}$, we have
$$
{\mathbb P}\{Y\cap S_i =\varnothing  \}\le 
(1-\vartheta)^M\le e^{-\vartheta M}.
$$
Hence,  the expected  number of sets $S_i\in{\mathcal S}$ disjoint from $Y_0$ is at most $Ne^{-\vartheta M}$. Choosing one additional point in each such set, we shall get a set $Y$ of cardinality $Ne^{-\vartheta M}+M$
intersecting all $S_i$. Puting $M=\lceil \vartheta^{-1}\log (N\vartheta)\rceil$, we get the desired bound.
\ep

Now, let $K\in\KR$.
Suppose that we can construct a probability measure $\mu$ on $\partial K$ such 
that for every $x\in \partial K$ and every $\varepsilon>0$, we have 
$\mu(S(x,\varepsilon))\ge p\varepsilon^{\frac{d-1}{2}}$ with some  $p>0$.

We take the Bronshtein--Ivanov net $X$ of $K$ constructed in Section~\ref{TBI}. 
Its  cardinality $N$ does not exceed $C(d)\ve^{-\frac{d-1}{2}}$, where $C(d)$ is of order $e^{O(d\log d)}$.
 Consider the  caps $S(x,\frac{\ve}{2}), x\in X$.
By Lemma~\ref{lem:transversal}, there exists a set $Y\subset\partial K$ of cardinality at most $\lceil 2^{\frac{d-1}{2}}p^{-1}\ve^{-\frac{d-1}{2}}\log (C(d) 2^{-\frac{d-1}{2}}p)\rceil +2^{\frac{d-1}{2}}p^{-1}\ve^{-\frac{d-1}{2}}$ that 
intersects each of those caps. If $p=e^{O(d)}$, then  the cardinality of $Y$  is of order $e^{O(d)}\ve^{-\frac{d-1}{2}}$.

\section{The construction of the measure}
\label{Tc}

Let $n$ be a positive integer (we shall need both $n=d$ and $n=d-1$).
Recall that for a convex body $K\subset {\mathbb R^n}$ containing the origin in 
its interior, its polar body $K^{\circ}\subset{\mathbb R^n}$ is defined by
$$
K^{\circ}=\{y\in{\mathbb R^n}:\langle x, y\rangle\le 1\quad\textrm{for 
all}\quad x\in K\}.
$$
We shall need the following well-known (but, in part, highly non-trivial) facts 
about the polar bodies:

\par\noindent\textit{Fact 1}.  If $K$ has a smooth boundary and is strictly convex, that is, $K$ contains no 
line segment on its boundary, then the 
relation $\langle x, x^*\rangle=1$, $x\in\partial K$, $x^*\in \partial 
K^{\circ}$, defines a continuous one to one mapping ${\,}^*$ from $\partial K$ 
to $\partial K^{\circ}$. The vector $x^*$ is just $\frac{\nu}{\langle x, 
\nu\rangle}$, where $\nu$ is the outer unit normal to $\partial K$ at $x$ (see 
\cite{Sch}, Corollary 1.7.3, page 40).

\par\noindent\textit{Fact 2}. For any convex body $K\subset{\mathbb R^n}$ 
containing the origin in its interior, we have 
$\textrm{vol}_n(K)\textrm{vol}_n(K^{\circ})\ge e^{O(n)}n^{-n}$ (see \cite{BM}, 
\cite{K}, \cite{NAZ}).

\par\noindent\textit{Fact 3}.  If $K$ is a   convex body  with the center of mass at the origin, then 
$$ 
\operatorname{vol}_n(K)\operatorname{vol}_n(K^{\circ})\le e^{O(n)}n^{-n}
$$
(see \cite{MP}).

\bl\label{dob}
Let    $K\subset {\mathbb R^d}$  be a convex body   containing  the origin in its interior 
and satisfying the Santal\'o bound  
$\operatorname{vol}_d(K)\operatorname{vol}_d(K^{\circ})\le e^{O(d)}d^{-d}$.
For any Borel set $S\subset \partial K$, define $S^*=\{x^*\in\partial 
K^{\circ}: x\in S\}$. Consider the ``cones" $C(S)=\{rx: x\in S, 0\le r\le 
1\}$ and $C(S^*)=
\{ry: y\in S^*, 0\le r\le 1\}$ and put
$$
\mu(S)=\frac{1}{2}\Big(\frac{\operatorname{vol}_d(C(S))}{\operatorname{vol}_d(K)}+\frac{
\operatorname{vol}_d(C(S^*))}{\operatorname{vol}_d(K^{\circ})}   \Big).
$$
Then $\mu$ is a probability  measure on $\partial K$ invariant under linear 
automorphisms of $\mathbb R^d$ and $\mu(S(x,\varepsilon))\ge 
e^{O(d)}\varepsilon^{\frac{d-1}{2}}$ for all $x\in\partial K$ and all $\ve\in(0,\frac{1}{2})$.
\el
\bp
The invariance of $\mu$ under linear automorphisms of ${\mathbb R^d}$ follows immediately from the general properties of the volume with respect to linear transformations and the relation 
$(TK)^{\circ}=(T^{-1})^*K^{\circ}$.

Apply an appropriate linear transformation to put the body $K$ in 
such a position that $x=x^*=e=(0,\dots,0,1)\in{\mathbb R^d}$. Then 
$S=S(x,\varepsilon)$ is given by $ \langle x, e\rangle\ge 1-\varepsilon$. Let 
$Q\subset e^{\perp} \cong {\mathbb R^{d-1}}$ be the convex body such that 
$(1-\varepsilon)e+Q$ is the cross-section of $K$ by the hyperplane $\{x: 
\langle x, e\rangle= 1-\varepsilon  \}$.
Let $\widetilde{K}=K\cap \{x: \langle x, e\rangle\le 1-\varepsilon  \}$.

\begin{figure}[ht]
\includegraphics[height=5cm]{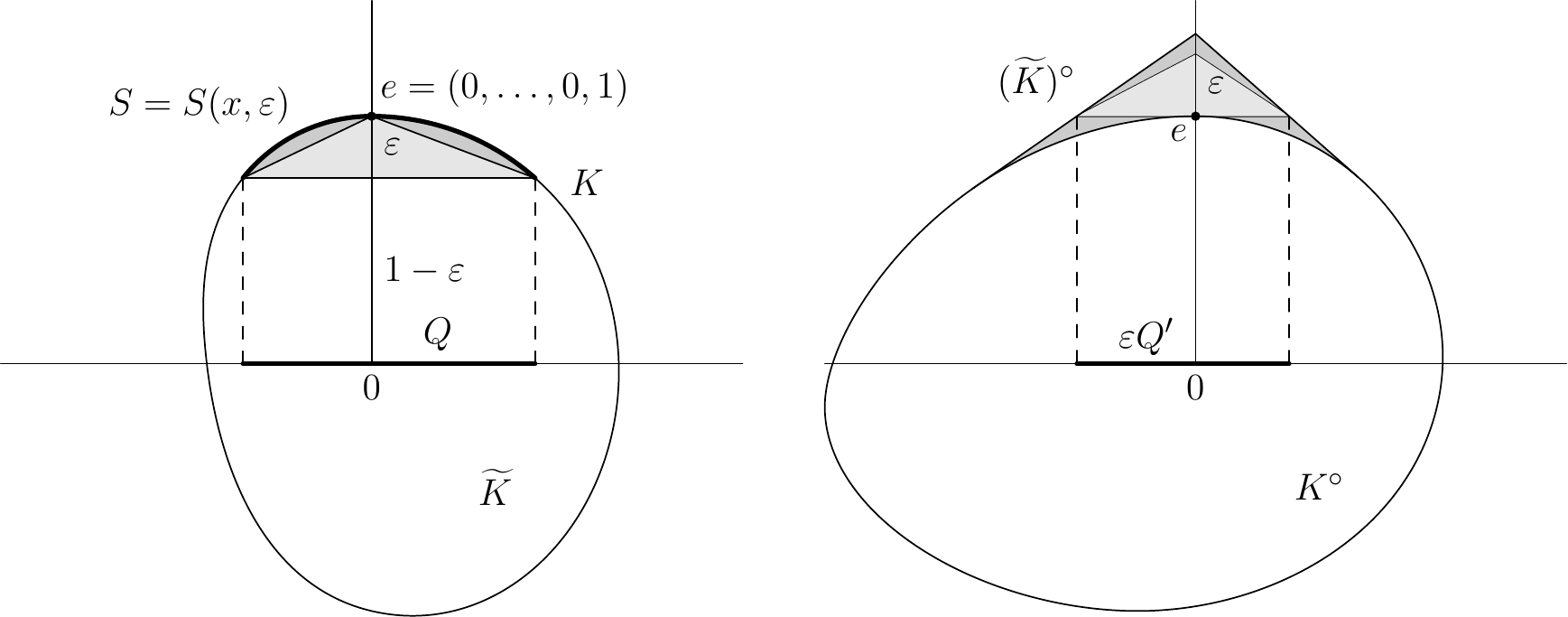}
\caption{The regions $K\setminus\widetilde K$ and $(\widetilde K)^\circ\setminus K^\circ$}
\label{fig4}
\end{figure}

Our first goal will be to show that 
$$
\operatorname{vol}_d(K\setminus \widetilde{K})\operatorname{vol}_d(
(\widetilde{K})^{\circ}\setminus K^{\circ})\ge
\frac{1}{d^2}\varepsilon^{d+1}\operatorname{vol}_{d-1}(Q)\operatorname{vol}_{d-1
}( Q'),
$$
where $Q'\subset e^\perp$ is the polar body to $Q$ in $\mathbb R^{d-1}$.

To this end,  note that $K\setminus \widetilde{K}$ contains the interior of the  pyramid
$\operatorname{conv}(\{e\}\cup (1-\varepsilon)e+Q)$ of height $\varepsilon$
with the base $(1-\varepsilon)e+Q$, so
$$
\operatorname{vol}_d(K\setminus \widetilde{K})\ge
\frac{1}{d}\varepsilon\operatorname{vol}_{d-1}(Q).
$$

We claim now that the interior of the pyramid
$\Pi=\operatorname{conv}\{(1+\varepsilon)e, e+\varepsilon Q'\}$ is
contained in $(\widetilde{K})^{\circ}\setminus  K^{\circ}$ (see Figure \ref{fig4}). Since
$ K^{\circ}\subset \{y:\langle y, e\rangle\le 1\}$, and
$\operatorname{int}\Pi\subset \{y: \langle y, e\rangle> 1\}$, it suffices 
to show
that $\Pi\subset (\widetilde{K})^{\circ}$.

To this end, take
$x\in \widetilde{K}$, and let $\langle x, e\rangle=
1-t\varepsilon$, $t\ge 1$, so $x=
(1-t\varepsilon)e+x'$, where $x'\in e^{\perp}$. 

\begin{figure}[ht]
\includegraphics[height=5cm]{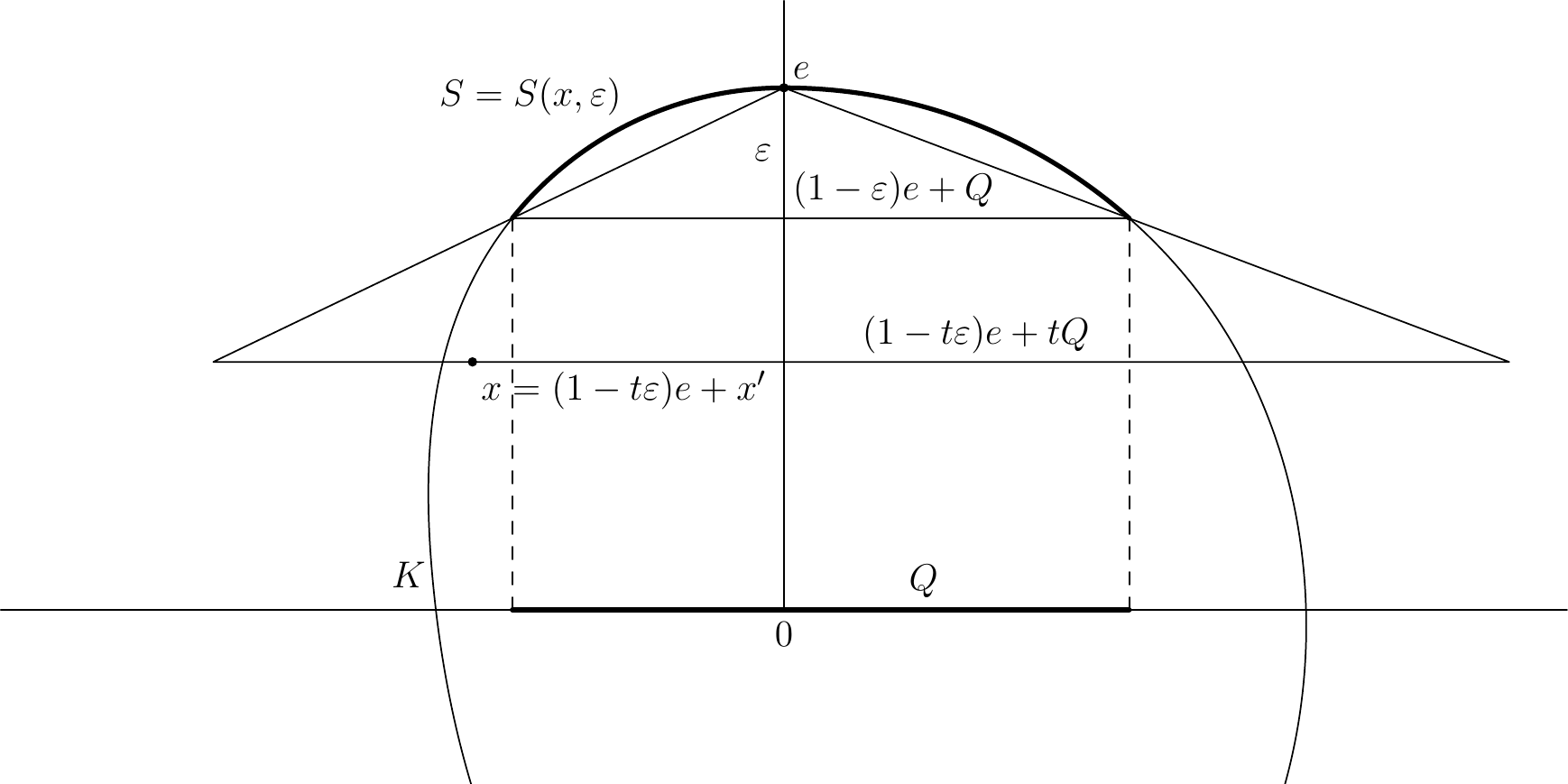}
\caption{The cross-section of $K$ by the hyperplane 
$\{x:\langle x,e\rangle=1-t\ve\}$ is contained in $tQ$}
\label{nefig5}
\end{figure}

Since $e\in K$, by the
convexity of $K$, $x'\in
tQ$ (see Figure \ref{nefig5}). Now, $\langle x, (1+\varepsilon)e\rangle =(1-t\ve)(1+\ve)\leq 1$, hence,
$(1+\varepsilon)e\in (\widetilde{K})^{\circ}$.
Let $y=e+\varepsilon y'$ with $y'\in Q'$. Then
$\langle x, y\rangle =1-t\varepsilon+\varepsilon \langle x',y'\rangle \leq 
1-t\ve+t\ve= 1$. Thus,
$e+\varepsilon Q'\subset (\widetilde{K})^{\circ}$. It follows by the convexity 
of $(\widetilde{K})^{\circ}$ that
$\Pi\subset (\widetilde{K})^{\circ}$, and, therefore, 
$$
\operatorname{vol}_d(
(\widetilde{K})^{\circ}\setminus K^{\circ})\ge\operatorname{vol}_d(
\Pi)=\frac{1}{d}\ve^d\operatorname{vol}_{d-1}(Q').
$$
Multiplying these two estimates, we get the desired inequality.

On the other hand, we have
$\operatorname{int}(K\setminus 
\widetilde{K})\subset C(S)\setminus (1-\varepsilon)C(S)$, and 
$\operatorname{int}((\widetilde{K})^{\circ}\setminus  K^{\circ})\subset 
(1-\varepsilon)^{-1}C(S^*)\setminus C(S^*)$. Hence, 
\begin{multline*}
\operatorname{vol}_d(K\setminus \widetilde{K})\operatorname{vol}_d( 
(\widetilde{K})^{\circ}\setminus K^{\circ})\le
\\
(1-(1-\varepsilon)^d)((1-\varepsilon)^{-d}-1)
\operatorname{vol}_d(C(S))\operatorname{vol}_d(C(S^*))\le
\\
e^{O(d)}\varepsilon^2\operatorname{vol}_d(C(S))\operatorname{vol}_d(C(S^*)).
\end{multline*}
Combining it with the previous estimate and using Fact 2, we get 
\begin{multline*}
\operatorname{vol}_d(C(S))\operatorname{vol}_d(C(S^*))\ge 
e^{O(d)}\varepsilon^{d-1}\operatorname{vol}_{d-1}(Q)\operatorname{vol}_{d-1}( 
Q')\ge
\\
e^{O(d)}\varepsilon^{d-1}(d-1)^{-(d-1)}.
\end{multline*}

Finally, since 
$\textrm{vol}_d(K)\textrm{vol}_d(K^{\circ})\le e^{O(d)}d^{-d}$, we get
\begin{multline*}
\mu(S)\ge
\frac{1}{2}\Big(\frac{\textrm{vol}_d(C(S))}{\textrm{vol}_d(K)}+\frac{
\textrm{vol}_d(C(S^*))}{\textrm{vol}_d(K^{\circ})}   \Big)\ge 
\\
\sqrt{\frac{\operatorname{vol}_d(C(S))\operatorname{vol}_d(C(S^*))}{
\textrm{vol}_d(K)\textrm{vol}_d(K^{\circ})}}\ge 
e^{O(d)}\varepsilon^{\frac{d-1}{2}},
\end{multline*}
as required.
\ep

This lemma, together with the discussion in Section  \ref{TR}, completes the proof of the theorem.

\end{document}